\documentclass[preprint,11pt]{elsarticle}

\usepackage{lipsum}
\makeatletter
\def\ps@pprintTitle{%
 \let\@oddhead\@empty
 \let\@evenhead\@empty
 \def\@oddfoot{}%
 \let\@evenfoot\@oddfoot}
\makeatother

\usepackage{amssymb}
\usepackage{amsthm}
\usepackage{amsmath}
\usepackage{algorithm}
\usepackage{algpseudocode}

\usepackage{stmaryrd}
\usepackage{tikz}

\newtheorem{theorem}{Theorem}
\biboptions{sort&compress}

\usepackage{xcolor}

\usepackage[caption=false]{subfig}
\usepackage{tikz}
\usepackage{verbatim}
\usepackage{color}
\usepackage{pgfplots}
\pgfplotsset{compat=newest}
\usetikzlibrary{positioning,shapes}
\usetikzlibrary{decorations.markings}
\usetikzlibrary{intersections}
\begin{document}

\begin{frontmatter}



\title{An algorithm based on continuation techniques for minimization problems with highly non-linear equality constraints}


\author[1]{E. Alberdi}
\address[1]{Matematika Aplikatua, Bilboko Ingeniaritza Eskola, University of the Basque Country UPV/EHU}

\author[2]{M. Anto\~{n}ana}
\address[2]{Matematika Aplikatua, Gipuzkoako Ingeniaritza Eskola, University of the Basque Country UPV/EHU}

\author[3]{J. Makazaga}
\address[3]{Konputazio Zientzia eta Adimen Artifiziala, Informatika Fakultatea, University of the Basque Country UPV/EHU}

\author[3]{A. Murua}
\address[3]{Konputazio Zientzia eta Adimen Artifiziala, Informatika Fakultatea, University of the Basque Country UPV/EHU}

\begin{abstract}
We present an algorithm based on continuation techniques that can be applied to solve numerically minimization problems with equality constraints.  We focus on problems with a great number of local minima which are hard to obtain by local minimization algorithms with random starting guesses. We are particularly interested in the computation of minimal norm solutions of underdetermined systems of polynomial equations. Such systems arise, for instance, in the context of the construction of high order optimized differential equation solvers.  By applying our algorithm, we are able to obtain 10th order  time-symmetric composition integrators with smaller 1-norm than any other integrator found in the literature up to now.

\end{abstract}

\end{frontmatter}



\section{Introduction}
We are concerned with the numerical solution of constrained minimization problems of the form
\begin{eqnarray}\label{eq:minimize}
&&\text{minimize } g(x)\\
\label{eq:constraints}
&&\text{subject to } f_j(x)=0, \quad j=1,\ldots,m,
\end{eqnarray}
where $g: \mathbb{R}^n \to \mathbb{R}$ is the objective function, and  $f_j: \mathbb{R}^n\to \mathbb{R}$, $j = 1,\ldots, m\leq n$,  are the equality constraints.  
We exploit the ability to follow implicitly defined curves to join a point that locally minimizes a problem with $k$ constraints with points that locally minimize the same problem but with an added constraint.
This way, the minimization problem with constraints can be afforded as a sequence of problems, starting from a problem with less constraints than the original, whose local minima can be known, and by adding step by step the rest of constraints.

A variety of local minimization methods for minimization problems with non-linear equality constrains are available. One of the ways consists in using a sequence of quadratic programming subproblems (SQP) which have to be solved successively. Since SQP methods were developed \cite{han1976, han1977, powell1978fast}, a great variety of research has been carried out focusing on this technique \cite{ izmailov, lawrence, liu2007, schitt1982, schitt1983}. Many SQP algorithms are combined with trust region methods. The origin of trust region methods can be found in the classical Levenberg-Marquardt method \cite{levenberg, morelevenberg}.  The first trust region strategy for equality constrained minimization problems is developed in \cite{avardi}. Different strategies to define trust region subproblems are defined in \cite{celis, Powell1990, yuan2015}. 

We focus on problems having a high number of local minima that are difficult to obtain by local minimization algorithms. 
We are interested in situations where the system  of equality constraints (\ref{eq:constraints}) is highly nonlinear, in the sense that linear  approximations of constraints are only valid in very small neighbourhoods. As a consequence, algorithms based on SQP methods and trust region algorithms may converge very slowly.

We are particularly interested in the computation of minimal norm solutions of under-determined systems of polynomial equations. With that aim, we consider constrained minimization problems of the form 
\begin{eqnarray}
&&\text{minimize } \|x\|^2\\
\label{eq:constraintspol0}
&&\text{subject to } p_j(x)=0, \quad j=1,\ldots,m \leqslant n,
\end{eqnarray}
where $\|x\|$ is the Euclidean norm of $x \in \mathbb{R}^n$ and each
$p_j(x)$ is a polynomial of degree $d_j$. More precisely,
we are interested in computing points that locally minimize the Euclidean norm $\|x\|$ under the polynomial constraints (\ref{eq:constraintspol0}).  We note that the problems that we have in mind are too complex to apply algorithms that guarantee the computation of all local minimizers of the Euclidean norm. Hence, we aim at computing as many local minimizers as possible, and select those points with smallest norm.
We stress that we do not particularly focus on the efficiency of the algorithms, but rather in the ability to compute (with reasonable computational resources) as many small norm solutions of (\ref{eq:constraintspol0}) as possible.

Coming back to the general problem (\ref{eq:minimize})--(\ref{eq:constraints}), 
one could try to solve that problem with a graduated optimization approach \cite{blake, chap1, chaudhuri}, by subsequently considering a sequence of minimization problems $P_{0}, P_{1}, \ldots,P_r$,  where for each $k=0,1,\ldots,r$ (for some $r< m$) the problem $P_k$ consists on finding $x_k \in \mathbb{R}^n$ that globally minimizes $g(x)$ subject to the constraints $x \in \mathcal{C}_k$, where 
 \begin{equation}
 \label{eq:Ck}
\mathcal{C}_k = \{x \in \mathbb{R}^n\ : \ f_j(x)=0 \mbox{ for } j=1,\ldots,m+k-r\}.
\end{equation}
Obviously, 
 \begin{equation*}
\min_{x \in \mathcal{C}_0} g(x)  \leq \min_{x \in \mathcal{C}_1} g(x)  \leq \cdots \leq \min_{x \in \mathcal{C}_r} g(x). 
\end{equation*}

 However, computing a guaranteed global minimum for each of the considered minimization problems $P_k$ may not be possible in practice. On the other hand,
 it may happen that, once a global minimizer $x_k^*$  of the problem $P_k$ is computed, no feasible points of the constraints $\mathcal{C}_k$ exists in the vicinity of $x_k^*$.
 
This motivates us to consider a graduated local minimization approach based on relating the solutions of a sequence $P_0,\ldots,P_r$ of local minimization  (rather than global minimization) problems. 
We propose an algorithm that starts by computing a set $\mathcal{S}_0$ of points that locally minimize $g(x)$ subject to $x \in \mathcal{C}_{0}$, and from them, by using path continuation techniques, subsequently computes, for each $k=1,\ldots,r$, a set $\mathcal{S}_k$ of points that locally minimize $g(x)$ subject to $x \in \mathcal{C}_{k}$.  The algorithm will be considered successful if the minimum of $g(x)$ over the points in $\mathcal{S}_{r}$ is close to the global minimum of the problem (\ref{eq:minimize})--(\ref{eq:constraints}).
We thus reduce the task of computing local minimizers of the original continuous constrained minimization problem to the task of exploring a discrete graph whose vertices are stationary points of the Lagrange functions of the constrained minimization problems. 

The plan of the paper is as follows.  The theoretical framework of the connecting algorithm is given in Section~\ref{s:2} and the precise algorithm proposed in the present work is described in Section~\ref{s:3}. 
In Section ~\ref{s:4}, we treat in detail the case where the objective function $g(x)$ is the Euclidean norm and the constraints are polynomials.  In Section~\ref{s:5}, we consider as benchmark problem the constrained minimization problems (corresponding to the construction of optimized 10th order composition integrators) considered in~\cite{1sofro_spaletta}.

\section{Theoretical framework of the algorithm}
\label{s:2}

In practice, we aim at computing a set of critical points of the function $g(x)$ on the smooth manifold defined by the constraints (\ref{eq:constraints}), or equivalently, 
a set of stationary points $(x,\lambda)$ (where $x \in \mathbb{R}^n$ and $\lambda=(\lambda_1,\ldots,\lambda_m) \in \mathbb{R}^m$) of the Lagrange function
\begin{equation}
\label{eq:L}
\mathcal{L}(x,\lambda) = g(x) + \sum_{j=1}^{m} \lambda_j\, f_j(x)
\end{equation}
such that the minimum of $g(x)$ over this set of stationary points $(x, \lambda)$ is expected to be close to the global minimum.  The main idea will be to begin by computing a set $\mathcal{S}_0$ of
stationary points of the Lagrange function 
\begin{equation*} g(x) + \sum_{j=1}^{m-r} \lambda_j\, f_j(x)
\end{equation*}
corresponding to the minimization problem $\mathcal{P}_0$, and from them, obtain a set $\mathcal{S}_1$ of
stationary points of the Lagrange function 
\begin{equation*} g(x) + \sum_{j=1}^{m-r+1} \lambda_j\, f_j(x)
\end{equation*}
corresponding to the minimization problem $\mathcal{P}_1$, and so on. 

More precisely, the algorithm proposed in the present subsequently obtains, for each $k=0,1,\ldots,r-1$, non-singular zeros\footnote{$z^*\in \mathbb{R}^{\ell}$ is said to be a non-singular zero of the map $F:\mathbb{R}^{\ell} \to \mathbb{R}^{\ell}$ if $F(z^*) = 0$ and the Jacobian matrix of $F(z)$ at $z=z^*$ is invertible}  of the map
\begin{equation}
\label{eq:Fk}
F_k(x,\lambda) = 
 \left(
\begin{matrix}
\nabla_x g(x)  + \sum_{j=1}^{m} \lambda_j \, \nabla_x f_j(x)\\
 f_1(x) \\
 \vdots\\
 f_{m-r+k}(x)\\
 \lambda_{m-r+k+1}\\
 \vdots \\
 \lambda_{m}
\end{matrix}
\right),
\end{equation}
with the aim of finally computing  non-singular zeros of the map
\begin{equation}
\label{eq:Fr}
F_r(x,\lambda) = \left(
\begin{matrix}
\nabla_x g(x)  + \sum_{j=1}^{m} \lambda_j \, \nabla_x f_j(x) \\
 f_1(x) \\
 \vdots\\
 f_{m}(x)
\end{matrix}
\right),
\end{equation}
(and hence stationary points of the Lagrange function (\ref{eq:L})).
Of course, among all the possible zeros of the map (\ref{eq:Fr}) that can be computed in that way, we are interested in those having small values of $g(x)$.

We will actually consider a more general class of problems consisting on the following: Consider a function $G:\mathbb{R}^{\ell} \to \mathbb{R}$ and a finite sequence $F_k:\mathbb{R}^{\ell} \to \mathbb{R}^{\ell}$, $k=0,1,\ldots,r$, of smooth maps with the property that $F_{k}(z)-F_{k-1}(z)$ identically vanish for $\ell-1$ components. We aim at computing the minimum of $G(z)$ over the set of non-singular zeros of $F_r(z)=0$.   In the particular case of the constrained minimization problem (\ref{eq:minimize})--(\ref{eq:constraints}), $\ell=n+m$, $z=(x,\lambda)$, $G(x,\lambda)=g(x)$, and the maps $F_k(x,\lambda)$ are given by (\ref{eq:Fk})--(\ref{eq:Fr}). 

We work under the assumption that the non-singular zeros of $F_0(z)$ are easier to compute than those of $F_r(z)$.
The main idea will be to compute, for $k=1,\ldots,r$ and for each non-singular zero $z_{k-1}$ of $F_{k-1}(z)$,  the non-singular zeros of  $F_k(z)$ that are connected (in a sense to be defined in Section~\ref{s:2}) to $z_{k-1}$.  This will allow us to reformulate the original problem in terms of an oriented graph whose vertices are non-singular zeros of the maps $F_k(z)$, $k=0,1,\ldots,r$.

Let us denote, for $k=0,1,\ldots,r$ the set of non-singular zeros of $F_k(z)$ as $\mathbb{V}_k$ .
 Recall also that we aim at computing a non-singular zero of $F_r(z)=0$ that minimizes $G(z)$.  
For each $k \in \{1,\ldots,r\}$, 
\begin{itemize}
\item Let $H_k:\mathbb{R}^{\ell} \to \mathbb{R}^{\ell-1}$ be such that the  components of $H_k(z)$ are the $\ell-1$ common components of $F_{k-1}(z)$ and $F_{k}(z)$.
\item Let $\mathcal{U}_k$ be the subset of regular points of 
\begin{equation}
\label{eq:Hk=0}
\{z \in \mathbb{R}^{\ell} \ : \ H_k(z)=0\}.
\end{equation}
(That is, the set of points $z^* \in \mathbb{R}^{\ell}$ satisfying that $H_k(z^*)=0$ and that the Jacobian matrix of $H_k(z)$ at $z=z^*$ has full rank.) For any point $z^* \in \mathcal{U}_k$, the implicit function theorem guarantees the existence of a unique smooth curve of points $z$ of $\mathcal{U}_k$ containing $z^*$. Clearly, any non-singular zero of either $F_{k-1}(z)$ or $F_k(z)$ belongs to $\mathcal{U}_k$.
\item We define an equivalence relation $\stackrel{k}{\sim}$ in $\mathbb{V}_{k-1} \cup \mathbb{V}_{k}$ as follows: $z' \stackrel{k}{\sim} z''$ if there is a curve in $\mathcal{U}_k$ containing both $z'$ and $z''$.
\end{itemize}

We will next consider an oriented graph $\mathbb{G}$ whose set of vertices is the union $\mathbb{V} = \cup_{k=0}^{r}\mathbb{V}_k$. 

The pair $(z',z'') \in \mathbb{V} \times  \mathbb{V}$ is a directed edge (oriented from $z'$ to $z''$) of $\mathbb{G}$, if for some $k \in \{1,\ldots,r\}$, $z' \in \mathbb{V}_{k-1}$, $z'' \in \mathbb{V}_{k}$ 
and $z' \stackrel{k}{\sim} z''$.

We will say that a zero $z_r$ of $F_r(z)$ can be reached from $z_0 \in \mathbb{V}_0$  if there is a path of directed edges connecting $z_0$ with $z_r$. That is, if there exists a sequence of points $z_k \in \mathbb{V}_k$, $k=1,\ldots,r-1$, such that, for  each $k \in \{1,\ldots,r\}$, $(z_{k-1}, z_{k})$ is an edge of $\mathbb{G}$. 


Obviously, there may be zeros of $F_r(z)$ that cannot be reached in that way from points in $\mathbb{V}_0$. One may hope that, with an appropriate choice of the sequence of maps $F_0,F_1,\ldots,F_r$, the minimum of  $G(z)$ over the zeros  that can be reached from $\mathbb{V}_0$ coincides with (or at least is close to) the minimum of  $G(z)$ over the set $\{z \in \mathbb{R}^{\ell}\ : F_r(z)=0\}$.

\section{Description of the proposed algorithm}
\label{s:3}

In what follows, we assume that each $\mathbb{V}_k$ is finite. This will be certainly the case if the components of each $F_{k}(z)$ are defined as polynomials.

In principle, one could compute all the vertices in $\mathbb{V}_0$ (all the non-singular zeros of $F_0(z)$), and successively compute, for $k=1,\ldots,r$, all the vertices of $\mathbb{V}_k$ that are $k$-equivalent to some vertex in $\mathbb{V}_{k-1}$. This would allows us to explore the whole oriented graph $\mathbb{G}$ and finally arrive to the minimum of $G(z)$ subject to $z \in \mathbb{V}_r$.

However, in the most problems of interest, and in particular, in the application examples considered in Section~\ref{s:5}, the number of elements of each $\mathbb{V}_k$ is exceedingly high for practical purposes. Due to that, we only explore $\mathbb{G}$ partially, by prescribing a threshold value $G_{\max}$, and actually constructing the subgraph of $\mathbb{G}$ corresponding to the subset of vertices $z$ satisfying that $G(z) \leq G_{\max}$. 

For each $k\in \{1,\ldots,r\}$, let $w_k:\mathbb{R}^{\ell} \to \mathbb{R}$ be the component of $F_k(z)$ that is missing from the $\ell-1$ components of $H_k(z)$.  In practice, the computation of the points in  $\mathbb{V}_{k}$ that are $k$-equivalent to a given point $z_{k-1} \in \mathbb{V}_{k-1}$  requires the ability to follow the (uniquely defined) curve of regular points of (\ref{eq:Hk=0}) containing $z_{k-1}$, while checking if $w_k(z)$ changes its sign along such a curve. This will guarantee that we compute all the non-singular zeros of $F_k(z)$ along that curve. (We may miss some singular zeros of $F_k(z)$, as $w_k(z)$ does not necessarily change its sign along the curve.)

\begin{figure}[tb]
\centering
\subfloat[Example of closed curve]{%
\resizebox*{6cm}{!}{
\begin{tikzpicture}[scale=0.8]
\draw[blue, ultra thick, name path=first] plot [smooth] coordinates {(-1,3) (-1.5,1) (3,2.5) (-3,-3)  };
 \draw[ultra thick,red, name path=third,
decoration={markings, mark=at position 0.08 with {\arrow[line width=1mm]{<}}},
decoration={markings, mark=at position 0.17 with {\arrow[line width=1mm]{<}}},
postaction={decorate}] plot [smooth cycle] coordinates {(-2,2) (1,1) (1.1,-0.8) (-5,-3) (-4,0)};
\fill [black] (1.1,-0.8) circle (0.1);
\node[right] at (1.3, -0.8) {$z_{k-1} \in \mathbb{V}_{k-1}$};
\node[right] at (0,-1.7) {\textcolor{red}{$\mathbf{H_k(z)=0}$}};
\node[left] at (-1.7,-1.5) {\textcolor{blue}{$\mathbf{w_k(z)=0}$}};
\draw [execute at begin node={\global\let\n=\n}, name intersections={of=first and third, total=\n}]  
      \foreach \i in {1,...,\n} {(intersection-\i) coordinate (red-blue-intersection-\i)};
    \foreach \i in {1,...,\n}
      \draw (red-blue-intersection-\i) circle [radius=3pt];
\end{tikzpicture}
}}
\hspace{5pt}
\subfloat[Curve limited by singular points]{%
\resizebox*{6cm}{!}{
\begin{tikzpicture}[scale=0.8]
\draw[blue, ultra thick, name path=first] plot [smooth] coordinates {(-1,3) (-1.5,1) (3,2.5) (-3,-3)  };
\draw [ultra thick,red, name path=second,
decoration={markings, mark=at position 0.47 with {\arrow[line width=1mm]{>}}},
decoration={markings, mark=at position 0.38 with {\arrow[line width=1mm]{<}}},
postaction={decorate}] (-2,2) to[out=45,in=115] (1,1) to[out=-180+115,in=10] (-5,-3);
\fill [black] (1.1,-0.3) circle (0.1);
\fill [purple](-2,2) circle (0.1);
\fill [purple](-5,-3) circle (0.1);
\node[right] at (1.3, -0.3) {$z_{k-1} \in \mathbb{V}_{k-1}$};
\node[right] at (0,-1.7) {\textcolor{red}{$\mathbf{H_k(z)=0}$}};
\node[left] at (-1.7,-1.5) {\textcolor{blue}{$\mathbf{w_k(z)=0}$}};
\draw [execute at begin node={\global\let\n=\n}, name intersections={of=first and second, total=\n}]  
      \foreach \i in {1,...,\n} {(intersection-\i) coordinate (red-blue-intersection-\i)};
    \foreach \i in {1,...,\n}
      \draw (red-blue-intersection-\i) circle [radius=3pt];
\end{tikzpicture}
}}
\caption{Illustration of typical steps of the algorithm. The intersection points are the elements in $\mathbb{V}_k$ that are $k$-equivalent to $z_{k-1}$} \label{fig:jarraipena}
\end{figure}
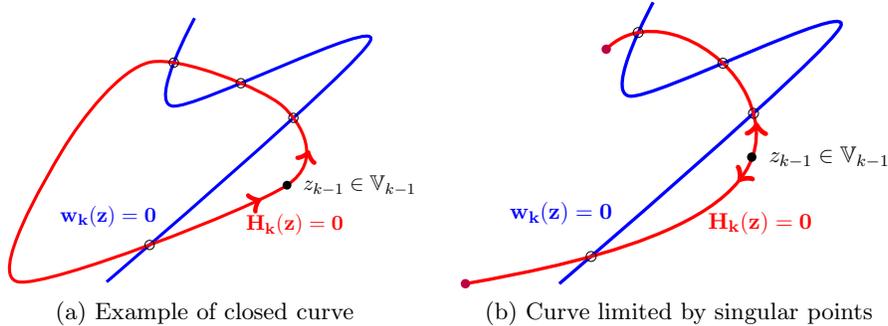

The algorithm proceeds in a step-by-step manner for $k=1,\ldots,r$ as follows: assume that a set $\mathcal{S}_{k-1}$ of non-singular zeros of $F_{k-1}(z)$ have been computed whose value of $G$ does not exceed $G_{\max}$.
\begin{itemize}
\item We begin by initializing $\mathcal{S}_{k}$ as the empty set.

\item  For each $z_{k-1} \in \mathcal{S}_{k-1}$, we find all the zeros of $F_k(z)$ that are $k$-equivalent to $z_{k-1}$ and include them  in $\mathcal{S}_{k}$.  In order to do that, the algorithm proceeds as follows. (See Figure~\ref{fig:jarraipena} for an illustration for the case where $\ell=2$.) Choose an orientation of the curve of regular points of (\ref{eq:Hk=0}) containing $z_{k-1}$, and follow that curve in the forward direction starting from $z_{k-1}$ while checking if  $w_k(z)$ changes its sign, until one or the following occurs: 
\begin{enumerate}
\item The curve arrives back to $z_{k-1}$. In that case, we stop searching for zeros of $F_k(z)$ that are $k$-equivalent to $z_{k-1}$.
\item A singular point of  (\ref{eq:Hk=0}) is reached. In that case, follow the curve in the backward direction starting from $z_{k-1}$  until some change of sign of $w_k(z)$ occurs or a singular point of  (\ref{eq:Hk=0}) is reached.
 \item A change of sign of $w_k(z)$ is found, which allow us to compute a zero $z_{k}$ of $w_k(z)$ along the curve, so that $F_k(z_k)=0$. In that case, if $z_k$ is already in $\mathcal{S}_k$,  then stop following that curve, as it would mean that such a curve has been already explored with a different starting point $z'_{k-1} \in \mathcal{S}_{k-1}$ (i.e., $z'_{k-1} \stackrel{k}{\sim} z_{k} \stackrel{k}{\sim} z_{k-1}$). Otherwise,  if $G(z_{k}) \leq G_{\max}$, then include $z_{k}$ in $\mathcal{S}_{k}$, and in any case, go on following the curve in the same direction, until another change of sign of $w_k(z)$ is found, or the curve arrives back to $z_{k-1}$, or a singular point of  (\ref{eq:Hk=0}) is reached.
\end{enumerate}

\end{itemize}

At the $r$th step, we get a set $\mathcal{S}_r$ on non-singular zeros of $F_r(z)$. We choose among them those points $z \in \mathcal{S}_r$ with smaller value $G(z)$. 

Some remarks on practical implementation aspects of our algorithm are made next:
\begin{itemize}
\item 
 We apply our own implementation (written in C) of a continuation algorithm that follows implicitly defined curves with a parametrization that is close to the arc length (referred to as pseudo-arc-length parametrization in~\cite{auto1, auto2, auto3}), which computes the zeros of a target function along the curve. 
\item It is considered that a singular point of (\ref{eq:Hk=0}) is reached if certain numerical difficulties are encountered to follow the curve further: The step-size used to advance along the curve is chosen in an adaptive way, by requiring that approximately the same number of simplified Newton iteration is needed at each step to follow the curve within a prescribed accuracy. A step is rejected (and a new step with smaller step size is tried next) if the iteration does not converge for a given step-size. Having too many consecutive rejected steps is considered as an indicator of the proximity of a singular point in the curve.
\item  In practice, we want to prevent wasting computing time by following a curve of excessive length with no zeros of the target function. We thus stop following a curve in a given direction (as in the case when a singular point of  (\ref{eq:Hk=0}) is achieved)
if the curve is followed for a (pseudo-) arc length larger than a prescribed 
 positive number $L_{\max}$.
 \item In order to avoid exploring the oriented graph $\mathbb{G}$ beyond practical computational limits, we consider (in addition to the use of the threshold $G_{\max}$) a maximum number $\ell_{\max}$ of zeros of $w_k(z)$ to be computed along each of the two orientations of the curve that connects a given point $z_{k-1} \in \mathcal{S}_{k-1}$ with points in $\mathcal{S}_{k}$.
\item The algorithm is easily implementable
in a network of processors and parallelizable in several ways: on the one hand,
the processors can treat separately each seed $z_{k-1}$ and on the other hand, each step $k$
may also start working as soon as the first element of $\mathcal{S}_k$ is obtained.
\end{itemize}

\section{Norm minimization subject to polynomial constraints}
\label{s:4}

In this section we describe the application of our method to the case where the constraints are polynomial and the objective function $g(x)$ is either the Euclidean norm. 

We thus consider constrained minimization problems of the form 
\begin{eqnarray}\label{eq:minimizepol}
&&\text{minimize } \sum_{i=1}^n \gamma_i^2\\
\label{eq:constraintspol}
&&\text{subject to } p_j(\gamma_1,\ldots,\gamma_{n})=0, \quad j=1,\ldots,m \leqslant n,
\end{eqnarray}
where 
$p_j(\gamma_1,\ldots,\gamma_{n})$ is a polynomial of degree $d_j$ in the variables $\gamma_1,\ldots,\gamma_{n}$.  

Assume that, for a prescribed positive integer $r$,  we are able to compute some subset $\mathcal{S}_{0}$ of the set of local minima of (\ref{eq:minimizepol}) subject to the constraints $p_j(\gamma_1,\ldots,\gamma_{n})=0$, $j=1,\ldots,m-r$.
Obviously, (\ref{eq:minimizepol})--(\ref{eq:constraintspol}) is of the form (\ref{eq:minimize})--(\ref{eq:constraints}), with  
$g(\gamma_1,\ldots,\gamma_n) = \|(\gamma_1, \ldots,\gamma_n)\|^2$, 
where $\|\cdot\|$ denotes the Euclidean norm. Hence,  the algorithm described in the previous section can be directly applied to that problem, giving as output a set $\mathcal{S}_{r}$ of stationary points of the Lagrangian function of the problem (\ref{eq:minimizepol})--(\ref{eq:constraintspol}).

However, we have found that it is advantageous rewriting the problem in an equivalent form by applying a technique that is standard in the numerical solution of polynomial system of equations. This consists in rewriting the constraints as a system of homogeneous polynomials (i.e., polynomials consisting on the sum of monomials of the same degree) with an additional indeterminate, say $\gamma_0$, and adding a constraint that confines the feasible points to a bounded set of $\mathbb{R}^{n+1}$.

 The problem  (\ref{eq:minimizepol})--(\ref{eq:constraintspol}) can be rephrased as the minimization of 
 \begin{equation*}
 \|(\gamma_1/\gamma_0, \ldots,\gamma_n/\gamma_0)\|^2 
\end{equation*}
  subject to the constrains
\begin{equation}
\label{eq:constraintsray}
p_j(\gamma_1/\gamma_0, \ldots, \gamma_n/\gamma_0)=0, \quad j=1,\ldots,m.\end{equation}
Clearly, if $(\gamma_0,\ldots,\gamma_n)$ is a solution of that constrained minimization problem, then  $(\mu\, \gamma_0,\ldots,\mu\, \gamma_n)$ is also a solution for each $\mu \in \mathbb{R}$. Hence, the same minimum is attained with the additional constraint $\|(\gamma_1, \ldots,\gamma_n)\|^2-R^2=0$
(with some fixed $R>0$), and in that case, the objective function can be replaced by $R^2/\gamma_0^2$. This is equivalent to the following:
\begin{eqnarray}\label{eq:minimizepolH}
&&\text{minimize } R^2/\gamma_0^2\\
\label{eq:constraintspolH}
&&\text{subject to }  P_j(\gamma_0, \gamma_1,\ldots,\gamma_{n})=0, \quad j=1,\ldots,m+1,
\end{eqnarray}
where 
\begin{equation}
\label{eq:P1}
P_1(\gamma_0,\gamma_1,\ldots,\gamma_n):=\sum_{i=1}^{n} \gamma_i^2-R^2, 
\end{equation}
and for each $j \in \{1,\ldots,m\}$,  we define the homogeneous polynomial
\begin{equation*}
P_{j+1}(\gamma_0, \ldots,\gamma_n)=\gamma_0^{d_j}p_j(\gamma_1/\gamma_0,\ldots,\gamma_n/\gamma_0).
\end{equation*}
where   $d_j$ is the degree of the polynomial $p_j(\gamma_1,\ldots,\gamma_n)$.

It is trivial to check that, provided that $(\gamma^*_0, \ldots,\gamma^*_n)$ is a solution of (\ref{eq:minimizepolH})--(\ref{eq:constraintspolH}), then, 
$(\gamma_1,\ldots,\gamma_n)=(\gamma^*_1/\gamma^*_0, \ldots, \gamma^*_n/\gamma^*_0)$
is a solution of (\ref{eq:minimizepol})--(\ref{eq:constraintspol}). This is also true if the objective function in (\ref{eq:minimizepolH}) is replaced by $-\gamma_0^2$.  

Thus, instead of computing a set of stationary points of the Lagrangian function of the problem (\ref{eq:minimizepol})--(\ref{eq:constraintspol}), we will compute a set of stationary points of the Lagrange function 
\begin{equation}
\label{eq:Lagrangian2}
 -\gamma_0^2 + \sum_{j=1}^{m+1} \lambda_j\, P_j(\gamma_0,\gamma_1,\ldots,\gamma_{n}).
\end{equation}
This has advantages from the point of view of numerical stability, and in addition, it allows us to move along curves that include points with vanishing $\gamma_0$ (which in the original formulation would correspond to points at infinity). 

Assume that, for a prescribed positive integer $r$,  we are able to compute some set $\mathcal{S}_{0}$ of points that locally minimize $-\gamma_0^2$ subject to the constraints $P_j(\gamma_0,\gamma_1,\ldots,\gamma_{n})=0$, $j=1,\ldots,m-r+1$. We will apply our algorithm with the sequence of polynomial maps  $F_k:\mathbb{R}^{\ell} \to \mathbb{R}^{\ell}$, $k=0,1,\ldots,r$ to be defined next in such a way that the points in $\mathcal{S}_{0}$ are zeros of $F_0(z)$, and the stationary points of (\ref{eq:Lagrangian2}) correspond to zeros of $F_r(z)$.

In order to avoid working with large values of the Lagrange multipliers, we replace the vector of Lagrange multipliers $(\lambda_1,\ldots,\lambda_{m+1})$ by a vector $(\lambda_0,\lambda_1,\ldots,\lambda_{m+1})$ such that 
$\sum_{j=0}^{m+1} \lambda_j^2 - R_{\lambda}^2$ with some fixed $R_{\lambda}>0$. The stationary points of the Lagrange function (\ref{eq:Lagrangian2}) are then obtained by solving $F_{r}(x,\lambda)=0$, where $x=(\gamma_0,\ldots,\gamma_n) \in \mathbb{R}^{n+1}$, $\lambda = (\lambda_0,\ldots,\lambda_{m+1}) \in \mathbb{R}^{m+2}$, and 
\begin{equation}
\label{eq:Fr2}
F_r(x,\lambda) = \left(
\begin{matrix}
\lambda_0^2 + \cdots + \lambda_{m+1}^2 - R_{\lambda}^2 \\
\lambda_0\,  \nabla_x g(x)  + \sum_{j=1}^{m+1} \lambda_j \, \nabla_x P_j(x) \\
 P_1(x) \\
 \vdots\\
 P_{m+1}(x)
\end{matrix}
\right),
\end{equation}
with $g(\gamma_0,\ldots,\gamma_n)=-\gamma_0^2$.
In addition, we define for $k=0,1,\ldots,r-1$,
\begin{equation}
\label{eq:Fk2}
F_k(x,\lambda) = \left(
\begin{matrix}
\lambda_0^2 + \cdots + \lambda_{m+1}^2 - R_{\lambda}^2 \\
\lambda_0\,  \nabla_x g(x)  + \sum_{j=1}^{m+1} \lambda_j \, \nabla_x P_j(x) \\
 P_1(x) \\
 \vdots\\
 P_{m-r+k+1}(x)\\
 \lambda_{m-r+k+2}(x)\\
 \vdots\\
 \lambda_{m+1}
\end{matrix}
\right).
\end{equation}
In the notation of Sections~2 and 3, 
we have that $\ell = n+m+3$, $z=(x,\lambda)$, $x = (\gamma_0,\ldots,\gamma_n)$, $\lambda = (\lambda_0,\ldots,\lambda_{m+1})$, and $G(x,\lambda) = -\gamma_0^2$.

Our algorithm will produce a set $\mathcal{S}_r$ of non-singular zeros 
\begin{equation*}
z^* =(\gamma_0^*,\ldots,\gamma_n^*,\lambda_0^*,\ldots,\lambda_{m+1}^*)
\end{equation*}
 of $F_r(z)$. Among them, those with largest value of $|\gamma_0^*|$,  correspond to candidates  $(\gamma_1,\ldots,\gamma_n):= (\gamma_1^*/\gamma_0^*,\ldots,\gamma_n^*/\gamma_0^*)$ to minimize (\ref{eq:minimizepol}) subject to the constraints (\ref{eq:constraintspol}).

\section{A benchmark problem: optimized 10th order time-symmetric composition methods}\label{s:5}

We are particularly interested in constrained optimization problems arising in the context of the construction of optimized differential equation solvers: Given an $n$-parameter family of integrators for some class of problems,  some requirements are imposed (for instance, that they attain certain order of convergence) that are equivalent to $m\leq n$ (typically polynomial) equations on the $n$ parameters of the family of integrators. If $m < n$, then an optimized integrator is chosen by  requiring to minimize some objective function, reflecting in some sense the quality of the integrators.  

\subsection{Statement of the  constrained minimization problem}

In the context of numerical integration of initial value problems of ordinary differential equations, composition methods refers to composing of a basic low order time-symmetric integrator~\cite{5mcLachlan,  12suzuki, 13yoshida} with different time-steps to obtain a higher order integration scheme. This technique is of particular interest in the context of geometric numerical integration.  One starts from a low accuracy basic integrator preserving some geometric features of the exact solution flow, and tries to increase the order of accuracy of the method while preserving some of the geometric properties of the basic method. We refer to~\cite{7geometric_hairer} and references therein for the interested reader. 

A $n$-stage composition method is determined by a vector $$x = (\gamma_1, \ldots,\gamma_n) \in \mathbb{R}^n$$
that has to satisfy certain polynomial equations (the so-called order conditions) for the method attaining a prescribed order of accuracy. A composition method is time-symmetric if 
\begin{equation}
\label{eq:sym}
\gamma_{j} - \gamma_{n-j+1} = 0, \quad \mbox{for} \quad j=1,\ldots,[n/2].
\end{equation}

In~\cite{1sofro_spaletta},  the construction  of time-symmetric composition integrators of order up to 10
(optimized in the sense of minimizing the 1-norm of $(\gamma_1, \ldots,\gamma_n)$) is considered. The most difficult cases treated in~\cite{1sofro_spaletta} correspond to 10th order integrators of different number $n$ of stages,  
$n=31$, $n=33$, and $n=35$, which have to satisfy, in addition to the symmetry conditions (\ref{eq:sym}),  $16$ polynomial equations in  the real variables $\gamma_1,\ldots,\gamma_n$.

In the present work, we consider, as a benchmark problem to test  our algorithm, the problem of determining the coefficients $\gamma_1,\ldots,\gamma_n$ of time-symmetric 10th order composition methods (obtained by composing a time-symmetric second order integrator) with minimal norm.  As in~\cite{1sofro_spaletta}, we consider three cases:  $n=31,33,35$,
subject to the $[n/2]$ symmetry conditions (\ref{eq:sym}) together with the 16 order conditions corresponding to time-symmetric 10th order composition methods~\cite{murua_sanzserna1999} (see also  \cite{7geometric_hairer}):
\begin{equation}\label{order2}
\sum_{k=1}^n \gamma_k - 1 = 0,
\end{equation}
\begin{equation}\label{order4}
\sum_{k=1}^n \gamma_k^{2i+1}=0, \quad i=1,2,3,4,
\end{equation}
\begin{equation}\label{order6}
 \sum_{k=1}^n \gamma_k^3\left(\sum_{l=1}^k{}^{'} \gamma_l\right)^2=0,
\end{equation}
\begin{equation}\label{order8}
\begin{aligned}
 \sum_{k=1}^n \gamma_k^5\left(\sum_{l=1}^k{}^{'} \gamma_l\right)^2=0 &,\\
\sum_{k=1}^n \gamma_k^3 \sum_{l=1}^k{}^{'} \gamma_l \sum_{m=1}^k{}^{'} \gamma_m^3  =0, \hspace{0.5cm}& \sum_{k=1}^n \gamma_k^3\left(\sum_{l=1}^k{}^{'} \gamma_l\right)^4=0,
\end{aligned}
\end{equation}

\begin{equation}\label{order10}
\begin{aligned}
\sum_{k=1}^n \gamma_k^7\left(\sum_{l=1}^k{}^{'} \gamma_l\right)^2=0 &,\\
\sum_{k=1}^n \gamma_k^5 \sum_{l=1}^k{}^{'} \gamma_l \sum_{m=1}^k{}^{'} \gamma_m^3  =0, \hspace{0.5cm}& \sum_{k=1}^n \gamma_k^3 \sum_{l=1}^k{}^{'} \gamma_l \sum_{m=1}^k{}^{'} \gamma_m^5  =0,\\
\sum_{k=1}^n \gamma_k^3 \left(\sum_{l=1}^k{}^{'} \gamma_l\right)^2 \sum_{m=1}^k{}^{'} \gamma_m^3 \sum_{n=1}^m{}^{'} \gamma_n  =0, \hspace{0.5cm}& \sum_{k=1}^n \gamma_k^5\left(\sum_{l=1}^k{}^{'} \gamma_l\right)^4=0,\\
\sum_{k=1}^n \gamma_k^3 \left(\sum_{l=1}^k{}^{'} \gamma_l\right)^3 \sum_{m=1}^k{}^{'} \gamma_m^3  =0, \hspace{0.5cm}& \sum_{k=1}^n \gamma_k^3\left(\sum_{l=1}^k{}^{'} \gamma_l\right)^6=0.\\
\end{aligned}
\end{equation}
\noindent Summation with a prime indicates that the last summation term is halved, that is, 
\begin{equation*}
\sum_{l=1}^{k}{}^{'} a_l = a_1+\cdots+a_{k-1} +\frac{a_k}{2}.
\end{equation*}

In~\cite{1sofro_spaletta}  they try to solve problem of determining the coefficients $\gamma_1,\ldots,\gamma_n$ of time-symmetric 10th order composition methods (obtained by composing a time-symmetric second order integrator) with minimal 1-norm. In order to avoid computational difficulties due to the non-smoothness of the 1-norm, We instead consider the minimization of the Euclidean norm of $x=(\gamma_1,\ldots,\gamma_n)$.

\subsection{Starting the algorithm}

The algorithm described in Section~\ref{s:3} requires as a previous step choosing a subset of the constraining equations (the {\em initial constraints}). In this sense, we choose the symmetry conditions (\ref{eq:sym}) and
the five simpler order conditions (\ref{order2})--(\ref{order4}).

As a first step, we need to compute a set of stationary points of the Lagrange function of the starting problem:
\begin{align}\label{eq:initialproblem1}
&\text{minimize } \sum_{j=1}^{n} \gamma_j^2\\
\label{eq:initialproblem2}
&\text{subject to (\ref{eq:sym})--(\ref{order4}).}
\end{align}

We begin by describing a procedure that allows us to compute a large set of  stationary points of the Lagrange function 
\begin{equation}
\label{eq:LagrangeF00}
 \sum_{j=1}^{n} \gamma_j^2 + \lambda_1 \left(\sum_{j=1}^{n} \gamma_j -1 \right) + 
 \sum_{k=1}^{4} \lambda_{k+1} \sum_{j=1}^{n} \gamma_j^{2k+1}.
\end{equation}

A procedure to compute stationary points of a similar Lagrange function 
is given in~\cite{mcLachlan2002}. The main idea is that the points satisfying (\ref{eq:sym})-(\ref{order2}) and having $n-5$ components equal to any of the other 5 components, are actually stationary points of (\ref{eq:LagrangeF00}).

One begins by choosing 5 positive integers $i_1 \geq i_2 \geq i_3 \geq i_4 \geq i_5$ such that $i_1+\cdots + i_5 = n$, set 
\begin{equation*}
\ell_1=0, \quad  \ell_{k+1}=\ell_{k}+i_{k}, \quad k=1,2,3,4,
\end{equation*}
and impose, in addition to  (\ref{eq:sym})-(\ref{order4}), that
\begin{equation}
\label{eq:equalities}
\gamma_{\ell_k+ j} =\gamma_{\ell_k+1}, \quad \mbox{for} \quad k=1,2,3,4,5, \quad 2 \leq j \leq i_k.
\end{equation}
This gives a system of $n$ equations for the $n$ unknowns $\gamma_1,\ldots,\gamma_n$. 

\begin{theorem}
\label{th:1}
Given 5 positive integers $i_1 \geq i_2 \geq i_3 \geq i_4 \geq i_5$ such that $i_1+\cdots + i_5 = n$, 
let $(\gamma_1^*, \ldots,\gamma_n^*) \in \mathbb{R}^n$ be obtained by permuting the components of a non-singular solution of the system of equations formed by (\ref{eq:sym})-(\ref{order2}) and (\ref{eq:equalities}).
Then, there exists a unique choice $(\lambda_1^*,\ldots,\lambda_5^*) \in \mathbb{R}^5$ of the Lagrange multipliers such that $(\gamma_1^*, \ldots,\gamma_n^*,\lambda_1^*,\ldots,\lambda_5^*) $ is a stationary point of (\ref{eq:LagrangeF00}).
\end{theorem}
Its proof is very similar to the proof of Proposition~2 in~\cite{mcLachlan2002} and will be omitted. 

Theorem~\ref{th:1}  can then be applied to obtain a large number of stationary points of the Lagrange function
\begin{equation}
\label{eq:LagrangeF0}
\begin{split}
 \sum_{j=1}^{n} \gamma_j^2 &+ \lambda_1 \left(\sum_{j=1}^{n} \gamma_j -1 \right) + 
 \sum_{k=1}^{4} \lambda_{k+1} \sum_{j=1}^{n} \gamma_j^{2k+1} \\
 &+ 
 \sum_{k=1}^{[n/2]} \lambda_{k+5} (\gamma_{m-k+1} - \gamma_{k})
 \end{split}
\end{equation}
of the constrained minimization problem  (\ref{eq:initialproblem1})--(\ref{eq:initialproblem2}). Indeed, if
\begin{equation*}
(\gamma_1, \ldots,\gamma_n,\lambda_1,\ldots,\lambda_5)  = (\gamma_1^*, \ldots,\gamma_n^*,\lambda_1^*,\ldots,\lambda_5^*) 
\end{equation*}
is a stationary point of (\ref{eq:LagrangeF00}) satisfying the symmetry conditions (\ref{eq:sym}), then 
\begin{equation*}
(\gamma_1, \ldots,\gamma_n,\lambda_1,\ldots, \lambda_5, \lambda_{5+1}, \ldots,\lambda_{5+[n/2]})  = (\gamma_1^*, \ldots,\gamma_n^*,\lambda_1^*,\ldots,\lambda_5^*, 0,\ldots,0) 
\end{equation*}
is a stationary point of the Lagrange function (\ref{eq:LagrangeF0}).

As an example, consider $n=31$,  and $(i_1,\ldots,i_5)=(16,8,3,2,2)$. There is only one real solution $(\gamma_1, \ldots,\gamma_n) \in \mathbb{R}^n$ of the system formed by (\ref{eq:simplereq}) and (\ref{eq:equalities}). That solution produces by permuting its components $31!/(16! 8! 3! 2! 2!) \approx 4*10^{15}$ different points, each of them giving rise to a different stationary point of (\ref{eq:LagrangeF00}). Among them, $15!/(8! 4!)=1351350$ points fulfill the symmetry conditions (\ref{eq:sym}), and each of them gives rise to a different stationary point of (\ref{eq:LagrangeF0}).  Observe that,  no point corresponding to a pattern 
$(i_1,\ldots,i_5)$ with more than one odd entry $i_k$ can satisfy the symmetry condition (\ref{eq:sym}).

For each such different pattern $(i_1,\ldots,i_5)$,  there is a large amount of stationary points of the Lagrange function (\ref{eq:LagrangeF0}) (obtained as solutions of the system of equations given by (\ref{eq:sym})--(\ref{order4}) and (\ref{eq:equalities})), and all of them could in principle be used as starting points for the algorithm proposed in Section~\ref{s:3}.
Since we aim at minimizing the 1-norm of $x=(\gamma_1,\ldots,\gamma_n)$ subject to some constraints, it makes sense to consider as starting points of our algorithm only those points that satisfy
\begin{equation}
\label{eq:1-normbound}
|\gamma_1| + \cdots + |\gamma_n| \leq N_{\max}  
\end{equation}
for an appropriately chosen threshold $N_{\max}$ for the 1-norm. In particular, we have considered $N_{\max}=7.5$.

However, the number of different starting points $x=(\gamma_1,\ldots,\gamma_n)$ that satisfy (\ref{eq:1-normbound}) with $N_{\max}=7.5$ (in addition to  (\ref{eq:sym})--(\ref{order4}) and (\ref{eq:equalities}) for some  pattern $(i_1,\ldots,i_5)$) is exceedingly high.

 The criteria we have adopted to reduce the number of initial seeds are motivated by the intended use of the solutions as coefficients of 10th order composition integrators: From one hand, it is an usual requirement for the coefficients of such integration methods~\cite{7geometric_hairer}
that 
\begin{equation}
\label{eq:condcumsum}
0 \leq  \sum_{k=1}^i \gamma_k \leq 1  \quad \mbox{for} \quad i=2,\ldots,n.
\end{equation}
  On the other hand, it can be seen that, in addition to minimizing the 1-norm of $x=(\gamma_1,\ldots,\gamma_n)$, it is desirable for a good integration method that the quantity
\begin{equation}
\label{eq:of2}
\max_{1 \leq i \leq [n/2]} |\gamma_{1} + \cdots + \gamma_{i-1} + \frac{1}{2}\gamma_{i}|  
\end{equation}
is as small as possible.

Among  all the solutions of (\ref{eq:sym})--(\ref{order4}) and (\ref{eq:equalities}) for some  pattern $(i_1,\ldots,i_5)$, we have chosen as starting set $\mathcal{S}_0$ of seeds for our algorithm those that simultaneously fulfill (\ref{eq:1-normbound}), (\ref{eq:condcumsum}), and the condition that (\ref{eq:of2}) is smaller than 0.8. That gives $1954677$ seeds for $n=31$, $4785415$ for $n=33$ and $5801580$ for $n=35$.  

\subsection{Minimizing the Euclidean norm}

Recall that, in each of the three cases  $n=31,33,35$, we aim at minimizing the Euclidean norm of $x=(\gamma_1,\ldots,\gamma_{n})$
subject to the $[n/2]$ symmetry conditions (\ref{eq:sym}) together with the 16 order conditions  (\ref{order2})--(\ref{order10}).
Clearly, this is a problem of the form (\ref{eq:minimizepol})--(\ref{eq:constraintspol}). In the previous section, we have obtained (in each of the cases $n=31,33,35$) a large set ${\mathcal S}_0$ of stationary points of of the Lagrange function corresponding to the constraints $P_j(x)=0$, $j=1,\ldots,m-r$, where $m=[n/2]+16$  and $r=11$. 
The homogeneous polynomials corresponding to that initial set of constraints is
\begin{equation}
\label{eq:simplereq}
\begin{aligned}
P_{j+1}(\gamma_0,\ldots,\gamma_n) &:= \gamma_{j} - \gamma_{n-j+1}, \quad \mbox{for} \quad j=1,\ldots,[n/2], \\
P_{[n/2]+1}(\gamma_0,\ldots,\gamma_n) &:= -\gamma_0 +\sum_{k=1}^{n} \gamma_k, \\
P_{[n/2]+i+1}(\gamma_0,\ldots,\gamma_n) &:= \sum_{k=1}^{n} \gamma_k^{2i+1}, \quad i=1,2,3,4.
\end{aligned}
\end{equation}
Observe that the remaining $r=11$ constraints (\ref{order6})--(\ref{order10}) are already in homogeneous form. 
As for the additional constraint (\ref{eq:P1}) introduced to work with homogeneous polynomials, we have chosen $R=4$ as the radius of the corresponding sphere. 
The objective function to be minimized under the constraints (\ref{eq:constraintspolH})--(\ref{eq:P1}) is $G(\gamma_0,\ldots,\gamma_n)=-\gamma_0^2$.

We have applied the algorithm described in Section \ref{s:3} with  $\ell_{\max}=10$ and $G_{\max}=1.4, 1.8, 1.9$ for $n=31,33,35$ respectively (see Section~3 for the meaning of the parameters $\ell_{\max}$ and $G_{\max}$), and after several weeks of computations in a cluster, we have obtained in each of the three cases a large number of points that locally minimize the Euclidean norm. 

\subsection{Final results}

In order to compare our results with those presented in~\cite{1sofro_spaletta}, 
 we have obtained local minimizers of the 1-norm that are close to the computed local minimizers of the Euclidean norm (in each of the cases $n=31,33,35$, we have considered twenty points with smaller 1-norm),  by application of Newton methods to the Lagrange formulation of the constrained minimization problem.
In all the cases, at least two significant digits of the components of $x$ remain unchanged, and the 1-norm is only slightly reduced. 
In each of the three cases $n=31,33,35$, we obtain solutions with smaller $1$-norm than the ones obtained in \cite{1sofro_spaletta} (see Table~\ref{table3_first_stage}).  The 1-norm of our best solutions and those in~\cite{1sofro_spaletta} are compared in Table~\ref{table2_1norm_comparison}.

\begin{table}
\begin{center}
\begin{tabular}{|c|c|c|c|}
\hline
	&  1-norm & Euclidean norm\\ \hline
	$n=31$ & 1 & 9 \\ \hline
	$n=33$ & 3 & 0 \\ \hline
	$n=35$ & 7 & 119\\
	\hline
\end{tabular}
\end{center}
\caption{Number of solutions with smaller 1-norm and Euclidean norm than those in \cite{1sofro_spaletta}.}\label{table3_first_stage}
\end{table}

\begin{table}
\begin{center}
\begin{tabular}{|c|c|c|c|}
\hline
	&1-norm of solutions in \cite{1sofro_spaletta}	&1-norm of our best solutions\\ \hline
	$n=31$ & $7.54471205180863$ & $7.386456254909627$ \\ \hline
	$n=33$ & $6.790022344309263$ & $6.680425940964748$ \\ \hline
	$n=35$ & $6.181326366704916$ & $5.863208397834587 $\\
	\hline
\end{tabular}
\end{center}
\caption{Comparison of the 1-norm between our best solutions and the ones obtained in \cite{1sofro_spaletta}.}\label{table2_1norm_comparison}
\end{table}

\section{Concluding remarks}\label{s:concluding}

We have presented an algorithm to solve numerically minimization problems with equality constraints having a great number of local minima that are difficult to obtain with local minimization algorithms. We put special emphasis in the computation of minimal norm solutions of under-determined systems of polynomial equations.
We have successfully tested the algorithm with a benchmark problem (corresponding to the construction of optimized integrators for ordinary differential equations) previously considered in~\cite{1sofro_spaletta}.

The key feature of the proposed algorithm is that an optimization problem in a continuous domain is reduced to an optimization problem
in a discrete graph. The algorithm heavily relies on the availability of a large number of local minima of a reduced problem (corresponding to a subset of the constraints). In the considered benchmark problem, this was possible thanks to the special structure of the reduced optimization problem. In situations where this is not the case, our algorithm could still be useful if local minima of a reduced problem could be computed  with local minimization algorithms much more easily than for the full problem.

\section*{Acknowledgments}

All authors have received funding from the European Union’s Horizon 2020 research and innovation program under the Marie Sklodowska-Curie grant agreement No 777778. All of them were also partially funded by the Basque Government Consolidated Research Group Grant IT649-13 on ``Mathematical Modeling, Simulation, and Industrial Applications (M2SI)'' and the Project of the Spanish Ministry of Economy and Competitiveness with reference MTM2016-76329-R (AEI/FEDER, EU).

This work has been possible thanks to the support of the computing infrastructure of the i2BASQUE academic network.

\bibliographystyle{elsarticle-num} 
\bibliography{paper_algorithm_arxiv}






\appendix

\section{Coefficients of the new methods}

In this section, we include the coefficients of the 10th order time-symmetric composition methods which have the minimum 1-norm obtained by our algorithm in each case: $n=31$, $n=33$ and $n=35$.


Our 10-order and $n=31$ stages method's coefficients:
\begin{equation}
\begin{array}{rrrrl} 
\gamma_1 &=& \gamma_{31} &=& \hskip 0.82em 0.112021591030629\\ 
\gamma_2 &=& \gamma_{30} &=& \hskip 0.82em 0.431725601490890\\ 
\gamma_3 &=& \gamma_{29} &=& -0.179522661652292\\ 
\gamma_4 &=& \gamma_{28} &=& \hskip 0.82em 0.120580123137540\\ 
\gamma_5 &=& \gamma_{27} &=&  -0.398625072360396\\ 
\gamma_6 &=& \gamma_{26} &=& \hskip 0.82em 0.178939708529781\\ 
\gamma_7 &=& \gamma_{25} &=& \hskip 0.82em 0.110380761851205\\ 
\gamma_8 &=& \gamma_{24} &=& \hskip 0.82em 0.122821075302122\\ 
\gamma_9 &=& \gamma_{23} &=& \hskip 0.82em 0.424853834201251\\ 
\gamma_{10} &=& \gamma_{22} &=& \hskip 0.82em 0.080402608153253\\ 
\gamma_{11} &=& \gamma_{21} &=& -0.152579616423119\\ 
\gamma_{12} &=& \gamma_{20} &=& -0.518863729554078\\ 
\gamma_{13} &=& \gamma_{19} &=& \hskip 0.82em 0.098430328190055\\ 
\gamma_{14} &=& \gamma_{18} &=& -0.347022983737523\\ 
\gamma_{15} &=& \gamma_{17} &=& \hskip 0.82em 0.144536650569654\\ 
 & & \gamma_{16} &=& \hskip 0.82em 0.543843562542057
\end{array}
\end{equation}

Our 10-order and $n=33$ stages method's coefficients:
\begin{equation}
\begin{array}{rrrrl} 
\gamma_1 &=& \gamma_{33} &=& \hskip 0.82em 0.099136878219969\\ 
\gamma_2 &=& \gamma_{32} &=& \hskip 0.82em 0.091805759677231\\ 
\gamma_3 &=& \gamma_{31} &=& \hskip 0.82em 0.459401983479601\\ 
\gamma_4 &=& \gamma_{30} &=& -0.020010940625404\\ 
\gamma_5 &=& \gamma_{29} &=& \hskip 0.82em 0.289568761201962\\ 
\gamma_6 &=& \gamma_{28} &=& \hskip 0.82em 0.037676477495504\\ 
\gamma_7 &=& \gamma_{27} &=&  -0.234223019629333\\ 
\gamma_8 &=& \gamma_{26} &=&  -0.531940341338964\\ 
\gamma_9 &=& \gamma_{25} &=& \hskip 0.82em 0.229077943954870\\ 
\gamma_{10} &=& \gamma_{24} &=& \hskip 0.82em 0.125254188184227\\ 
\gamma_{11} &=& \gamma_{23} &=& \hskip 0.82em 0.154215725364726\\ 
\gamma_{12} &=& \gamma_{22} &=& \hskip 0.82em 0.095409688982420\\ 
\gamma_{13} &=& \gamma_{21} &=& \hskip 0.82em 0.048476867552146\\ 
\gamma_{14} &=& \gamma_{20} &=& -0.296771552754660\\ 
\gamma_{15} &=& \gamma_{19} &=& -0.337160630892827\\ 
\gamma_{16} &=& \gamma_{18} &=& \hskip 0.82em 0.011840660098572\\ 
 & & \gamma_{17} &=& \hskip 0.82em 0.556483102059918
\end{array}
\end{equation}

Our 10-order and $n=35$ stages method's coefficients:
\begin{equation}
\begin{array}{rrrrl} 
\gamma_1 &=& \gamma_{35} &=& \hskip 0.82em 0.100117054165055\\ 
\gamma_2 &=& \gamma_{34} &=& \hskip 0.82em 0.159849233601330\\ 
\gamma_3 &=& \gamma_{33} &=& \hskip 0.82em 0.316881415877955\\ 
\gamma_4 &=& \gamma_{32} &=& -0.221896402036101\\ 
\gamma_5 &=& \gamma_{31} &=& -0.231034183177538\\ 
\gamma_6 &=& \gamma_{30} &=& \hskip 0.82em 0.076265548489175\\ 
\gamma_7 &=& \gamma_{29} &=& \hskip 0.82em 0.110652300072783\\ 
\gamma_8 &=& \gamma_{28} &=& \hskip 0.82em 0.129556002817133\\ 
\gamma_9 &=& \gamma_{27} &=& \hskip 0.82em 0.094866828518147\\ 
\gamma_{10} &=& \gamma_{26} &=& \hskip 0.82em 0.114094318414488\\ 
\gamma_{11} &=& \gamma_{25} &=& \hskip 0.82em 0.255254772501037\\ 
\gamma_{12} &=& \gamma_{24} &=& \hskip 0.82em 0.070625655529692\\ 
\gamma_{13} &=& \gamma_{23} &=& -0.176094652551014\\ 
\gamma_{14} &=& \gamma_{22} &=& \hskip 0.82em 0.041045831082866\\ 
\gamma_{15} &=& \gamma_{21} &=& -0.210904961303419\\ 
\gamma_{16} &=& \gamma_{20} &=& -0.375871900390575\\ 
\gamma_{17} &=& \gamma_{19} &=& \hskip 0.82em 0.049098633077334\\ 
 & & \gamma_{18} &=& \hskip 0.82em 0.394989010623301
\end{array}
\end{equation}

\end{document}